\documentclass{amsart}
\usepackage{amssymb}
 
\newtheorem{theorem}{Theorem} 
\newtheorem{lemma}[theorem]{Lemma}

\theoremstyle{definition}

\theoremstyle{remark}

\newcommand{\rank}{\textit{\rm rank}}

\renewcommand{\>}{\rangle}

\DeclareSymbolFont{AMSb}{U}{msb}{m}{n}
\DeclareMathSymbol{\F}{\mathbin}{AMSb}{"46}
\DeclareMathSymbol{\N}{\mathbin}{AMSb}{"4E}
\DeclareMathSymbol{\Z}{\mathbin}{AMSb}{"5A}
\DeclareMathSymbol{\R}{\mathbin}{AMSb}{"52}
\DeclareMathSymbol{\C}{\mathbin}{AMSb}{"43}

\begin{document} \title[Minimal Generators for Symmetric Ideals]
{Minimal Generators for Symmetric Ideals}

 \author{Christopher J. Hillar}
\address{Department of Mathematics, Texas A\&M University, College Station, TX 77843.}
\email{chillar@math.tamu.edu}

 \author{Troels Windfeldt}
\address{Department of Mathematics, University of Copenhagen, Denmark.}
\email{windfeldt@math.ku.dk}

 \thanks{The work of the first author is supported under an NSF Graduate Research Fellowship.}






\mbox{\ } \vspace{-7pt}

\maketitle

Let $R = K[X]$ be the polynomial ring in infinitely many indeterminates $X$ 
over a field $K$. Write ${\mathfrak S}_{X}$ (resp. ${\mathfrak S}_{N})$
for the symmetric group of $X$ (resp. $\{1,\ldots,N\}$) and $R[{\mathfrak S}_{X}]$
for its (left) group ring, which acts naturally on $R$.  An \textit{invariant ideal} 
$I \subseteq R$ is an $R[{\mathfrak S}_{X}]$-submodule of $R$.  
Aschenbrenner and Hillar recently proved \cite{AH} that all invariant ideals
are finitely generated over $R[{\mathfrak S}_{X}]$.  
They were motivated by finiteness questions in
chemistry \cite{AH} and algebraic statistics \cite{SturmSull}.

In proving the Noetherianity of $R$, it was shown that 
an invariant ideal $I$ has a special, finite set of generators
called a \textit{minimal Gr\"obner basis}.  However, the more basic question
whether $I$ is always cyclic (already asked by J. Schicho) was left unanswered in \cite{AH}.  
Our result addresses a generalization of this important issue.



\begin{theorem}\label{mainthm}
For every positive integer $n$, there are invariant ideals of $R$ generated by $n$ polynomials 
which cannot have fewer than $n$ $R[{\mathfrak S}_X]$-generators.
\end{theorem}

In what follows, we work with the set
$X = \{x_1,x_2,x_3,\ldots\}$, although as remarked 
in \cite{AH}, this is not really a restriction.  In this case,
${\mathfrak S}_X$ is naturally identified with ${\mathfrak S}_\infty$, the
permutations of the positive integers, and $\sigma x_i = x_{\sigma i}$ for
$\sigma \in {\mathfrak S}_\infty$.


Let $M$ be a multiset of positive integers and let $i_1,\ldots,i_k$ be the list of its distinct elements,
arranged so that $m(i_1) \geq \cdots \geq m(i_k)$, where 
$m(i_j)$ is the multiplicity of $i_j$ in $M$.  
The \textit{type} of $M$ is the vector $\lambda(M) = (m(i_1),m(i_2),\ldots,m(i_k))$.
For instance, the multiset $M = \{1,1,1,2,3,3\}$ has type $\lambda(M) = (3,2,1)$.
Multisets are in bijection with monomials of $K[X]$.  Given $M$, 
we can construct the monomial:
\begin{equation*}
\mathbf{x}_{M}^{\lambda(M)} = \prod_{j =1}^k x_{i_j}^{m(i_j)}.
\end{equation*}
Conversely, given a monomial, the 
associated multiset is the set of indices appearing in it, along with
multiplicities.  The action of ${\mathfrak S}_\infty$ on monomials coincides with
the natural action of ${\mathfrak S}_\infty$ on multisets $M$ and this action 
preserves the type of a multiset (resp. monomial).
We also note the following elementary fact.
\begin{lemma}\label{infpermfiniteperm}
Let $\sigma \in {\mathfrak S}_\infty$ and $f \in R$.  Then there
exists a positive integer $N$ and $\tau \in {\mathfrak S}_N$ such that 
$\tau f = \sigma f$.
\end{lemma}


Theorem \ref{mainthm} is a direct corollary of the following result.

\begin{theorem}\label{genthm}
Let $G = \{g_1,\ldots,g_n\}$ be a set of monomials of degree $d$ 
with distinct types and fix a matrix $C = (c_{ij}) \in K^{n \times n}$ of rank $r$.  Then the
submodule $I = \<f_1,\ldots,f_n\>_{R[{\mathfrak S}_\infty]} \subseteq R$ generated by the $n$ 
polynomials, $f_{j} = \sum_{i = 1}^{n} c_{ij} g_i, (j = 1,\ldots, n),$
cannot be generated with fewer than $r$ polynomials.
\end{theorem}

\begin{proof}
Suppose that $p_1,\ldots,p_k$ are generators for $I$; 
we prove that $k \geq r$.  
Since each $p_l \in I$, it follows that each is a linear combination, 
over $R[{\mathfrak S}_\infty]$, of monomials in $G$.  Therefore, each monomial 
occurring in $p_l$ has degree at least $d$, and, moreover, any degree $d$ monomial
in $p_l$ has the same type as one of the monomials in $G$.

Write each of the monomials in $G$ in the form $g_i = \mathbf{x}_{M_i}^{\lambda_i}$ for
multisets $M_1,\ldots,M_n$ with corresponding distinct types $\lambda_1,\ldots,\lambda_n$,
and express each generator $p_l$ as:
\begin{equation}\label{pieqs}
p_l = \sum_{i=1}^{n} \sum_{\lambda(M) = \lambda_i} u_{i l M} \mathbf{x}_{M}^{\lambda_i} \  +  \ q_l, 
\end{equation}
in which $u_{i l M} \in K$ with only finitely many of them nonzero, 
each monomial in $q_l$ has degree larger than $d$,
and the inner sum is over multisets $M$ with type $\lambda_i$.

Since each polynomial in $\{f_1,\ldots,f_n\}$ is a finite linear combination of the $p_l$, and
since only finitely many integers are indices of monomials 
appearing in $p_1,\ldots,p_k$, we may pick $N$ large enough so that all of 
these linear combinations can be expressed with coefficients in the subring $R[{\mathfrak S}_N]$
(c.f. Lemma \ref{infpermfiniteperm}).
Therefore, we have,
\begin{equation}\label{sumpi}
\begin{split}
f_j = \ &  \sum_{l=1}^{k} \sum_{\sigma \in {\mathfrak S}_N} s_{l j \sigma} \sigma p_l, \\
\end{split}
\end{equation}
for some polynomials $s_{l j \sigma} \in R$.  Substituting (\ref{pieqs}) into (\ref{sumpi}) gives 
us that 
\begin{equation*}
\begin{split}
f_j  = \  &  \sum_{l=1}^{k} \sum_{\sigma \in {\mathfrak S}_N} 
\sum_{i=1}^{n} \sum_{\lambda(M) = \lambda_i} 
v_{l j \sigma}  u_{i l M} \mathbf{x}_{\sigma M}^{\lambda_i} \  + h_j, \\
\end{split}
\end{equation*}
in which each monomial appearing in $h_j \in R$ has degree greater than $d$ and $v_{l j \sigma}$
is the constant term of $s_{l j \sigma}$.  Since each $f_j$ has degree $d$, we 
have that $h_j = 0$.  Thus, 
\begin{equation*}
\begin{split}
\sum_{i = 1}^{n} {c_{ij} \mathbf{x}_{M_i}^{\lambda_i}} = \ &  \sum_{l=1}^{k} \sum_{\sigma \in {\mathfrak S}_N} 
\sum_{i=1}^{n} \sum_{\lambda(M) = \lambda_i} 
v_{l j \sigma}  u_{i l M} \mathbf{x}_{\sigma M}^{\lambda_i}. \\
\end{split}
\end{equation*}
Next, for a fixed $i$,
take the sum on each side in this last equation of the coefficients of monomials with 
the type $\lambda_i$.  This produces the $n^2$ equations:
\begin{equation*}\label{sumeq}
c_{i j} =   \sum_{l=1}^{k} \sum_{\sigma \in {\mathfrak S}_N} 
 \sum_{\lambda(M) = \lambda_i} 
v_{l j \sigma}  u_{i l M} \\
=  \sum_{l=1}^{k} \left( \sum_{\lambda(M) = \lambda_i}  u_{i l M} \right)\left(\sum_{\sigma \in {\mathfrak S}_N} 
v_{l j \sigma} \right) \\
=   \sum_{l=1}^{k} U_{i l} V_{l j},
\end{equation*}
in which $U_{i l} = \sum_{\lambda(M) = \lambda_i}  u_{i l M}$ and $V_{l j} = \sum_{\sigma \in {\mathfrak S}_N}   v_{l j \sigma}$.  Set $U$ to be the $n \times k$ matrix $(U_{i l})$ and similarly let 
$V$ denote the $k \times n$ matrix $(V_{l j})$.
These  $n^2$ equations are represented by the equation 
$C = UV$, leading to the following chain of inequalities:
\[ r = \rank(C) = \rank(UV) \leq \min \{\rank(U), \rank (V)\} \leq \min \{n,k\} \leq k.\]
Therefore, we have $k \geq r$, and this completes the proof.
\end{proof}



\begin{thebibliography}{99}

\bibitem{AH}
M. Aschenbrenner and C. Hillar, \emph{Finite generation of symmetric ideals}, 
Trans. Amer. Math. Soc., to appear.

\bibitem{SturmSull}
B. Sturmfels and S. Sullivant, \emph{Algebraic factor analysis: tetrads, 
pentads and beyond}, preprint. (math.ST/0509390).


\end{thebibliography}
\end{document}